\documentclass[journal]{IEEEtran}
\usepackage{amssymb}
\ifCLASSINFOpdf
\else
\usepackage[dvips]{graphicx}
\fi
\usepackage{graphics}
\usepackage{subfigure}
\usepackage{graphicx}
\usepackage{algorithm}
\usepackage{algorithmic}
\usepackage{graphics}
\usepackage{subfig}
\usepackage{balance}
\usepackage{amsmath}

\newcommand{\EE}{\mathbb{E}}
\newtheorem{theorem}{Theorem}
\newtheorem{definition}{Definition}
\newtheorem{lemma}{Lemma}

\begin{document}
\title{Non-ergodic Complexity of Convex  Proximal Inertial   Gradient Descents}
\author{Tao Sun, Linbo Qiao, Dongsheng Li\thanks{This work was
supported by the Major State Research Development Program (2016YFB0201305).

Part of a short version of the paper was accepted as a part of a regular  conference  proceeding  paper in   AAAI 2019.

The authors are  with
College of Computer, National University of Defense Technology,
Changsha, 410073, Hunan,  China (e-mails: \texttt{nudtsuntao@163.com}, \texttt{qiao.linbo@nudt.edu.cn},  \texttt{dsli@nudt.edu.cn}).

Dongsheng Li is the corresponding author.
}}

\maketitle

\begin{abstract}
The proximal inertial  gradient descent   is    efficient  for the composite minimization and
applicable for broad of machine learning problems.  In this paper, we revisit the computational complexity of  this algorithm  and  present other novel results, especially on the convergence rates of the objective  function values. The non-ergodic $O(1/k)$ rate is proved for  proximal  inertial gradient descent with  constant stepzise when the objective function is coercive. When the  objective function fails to promise coercivity, we prove the sublinear rate with diminishing inertial parameters. In the case that the objective function satisfies optimal strong convexity condition (which is much weaker than the strong convexity), the linear convergence is proved with much larger  and general stepsize than previous literature. We also extend our results to the multi-block version and present the computational complexity. Both cyclic and stochastic index selection strategies are considered.
\end{abstract}

\begin{IEEEkeywords}
Convex proximal inertial gradient  descent; Heavy-ball method;  Non-ergodic convergence rates; Block coordinate descent; Computational complexity
\end{IEEEkeywords}

\IEEEpeerreviewmaketitle

\section{Introduction}
This paper is devoted to study following composite minimization
\begin{equation}\label{model}
    \min_{x} \left\{F(x):=f(x)+g(x)\right\},
\end{equation}
where $f$ is differentiable, and $\nabla f$ is Lipschitz continuous with $L$, and $g$ is proximable. The  \b{p}roximal \b{i}nertial \b{g}radient \b{d}escent (PIGD) for   problem \eqref{model} performs the following iteration
 \begin{align}\label{scheme}
     x^{k+1}=\textbf{prox}_{\gamma_k g}[x^k-\gamma_k \nabla f(x^k)+\beta_k(x^k-x^{k-1})],
\end{align}
where $\gamma_k$ is the stepsize and $\beta_k$ is the inertial parameter.   The PIGD is proposed  in \cite{ochs2014ipiano} for nonconvex optimization and called as iPiano (\b{i}nertial \b{p}rox\b{i}mal \b{a}lgorithm for \b{n}onconvex \b{o}ptimization) in that paper.
In this paper, we consider the convex scenario and thus we use the name PIGD rather than iPiano.
PIGD is closely related to two classical algorithms: the forward-backward splitting method \cite{combettes2005signal} (when $\beta_k\equiv 0$ in \eqref{scheme}) and heavy-ball method \cite{polyak1964some} (when $g\equiv 0$ in \eqref{scheme}). PIGD is a combination of   forward-backward splitting method and  heavy-ball method. However, different from forward-backward splitting, the sequence  generated by PIGD is not Fej\'{e}r monotone due to the inertial term $\beta_k(x^k-x^{k-1})$.  This brings troubles in proving the convergence rates for the objective function values in the convex case. Noting that the heavy-ball method is a special form of PIGD, the difficulty also exists in analyzing the complexity of heavy-ball method. In the existing literatures, the sublinear convergence rate of the heavy-ball was established only in the sense of ergodicity.  In this paper, we propose a novel Lyapunov function to address this issue, and prove the non-ergodic convergence rates.
\subsection{Interpretation by Dynamical Systems}
Recently, there has been increasing interests in using ODEs to model and understand the first-order iterative optimization algorithms \cite{su2014differential,sun2017asynchronous,wibisono2016variational,attouch2018fast}.  This is because the ODE  is heavily related to numerical optimization; if  step sizes are  small enough so that the trajectory converges to a curve given by an ODE. The  well-established theory of
ODEs always  offers deeper insights for optimization algorithms, which has led to various novel and  interesting findings. More importantly, ODEs motivate  us to construct the proper Lyapunov function for the optimization algorithms. In fact, the  dynamical system for modeling  the heavy-ball method ($g\equiv 0$  in \eqref{scheme}) has been given in \cite{alvarez2000minimizing}:
\begin{equation}\label{ODE1}
    \ddot{x}(t)+\alpha \dot{x}(t)+\nabla f(x(t))=\textbf{0}
\end{equation}
for some $\alpha>0$. However, for PIGD, the corresponding  dynamical system has still  been unknown; we remedy this case. To introduce the dynamics, we need a definition given in \cite{su2014differential}. The the proximal subgradient is defined as
\begin{align*}
G_s(x) \triangleq \frac{x - \mbox{argmin}_{z}\left( \|z - (x-s\nabla f(x))\|^2/(2s) + g(z) \right) }{s}.
\end{align*}
\begin{definition}
A Borel measurable function $G(x, p;F)$ defined on $\mathbb{R}^n\times \mathbb{R}^{n}$ is said to be a directional subgradient of $F$ if
\begin{align*}
G(x, p) \in \partial F(x), ~~~ \langle G(x, p), p\rangle = \sup_{\xi\in \partial F(x)}\langle \xi, p \rangle
\end{align*}
for all $x, p$. Here, $G(x):=\lim_{s\rightarrow 0}G_{s}(x)$.
\end{definition}

The authors in \cite{su2014differential} pointed out that the existence of a directional subgradient can be guaranteed by  a lemma given in \cite{rockafellar2015convex}. If $g$ is differentiable, it is easy to verify that
\begin{align}\label{reason}
\langle G(x(t),\dot{x}(t)),\dot{x}\rangle=\frac{dF(x(t))}{dt}=\langle \nabla F(x(t)),\dot{x}(t)\rangle.
\end{align}
Equation \eqref{reason} indicates that $G(x(t),\dot{x}(t))$ can present  some  property of $\nabla F(x(t))$.
Noticing that $G(x(t),\dot{x}(t))$ exists in the nonsmooth case,  we thus consider replacing $\nabla f(x(t))$  with  $G(x(t),\dot{x}(t))$ in \eqref{ODE1} and propose the following system
\begin{equation}\label{ODE-PIAG}
    \ddot{x}(t)+\alpha \dot{x}(t)+G(x(t),\dot{x}(t))=\textbf{0},~~\,~\alpha>0.
\end{equation}
A similar dynamical system  related to the FISTA \cite{beck2009fast} and Nesterov's acceleration gradient method \cite{nesterov2013introductory} given in \cite{su2014differential} presented as
\begin{equation} \label{ODE-Nes}
    \ddot{x}(t)+\frac{3}{t} \dot{x}(t)+G(x(t),\dot{x}(t))=\textbf{0}.
\end{equation}
The difference between \eqref{ODE-PIAG} and  \eqref{ODE-Nes} lies in the coefficient for $\dot{x}(t)$: one is constant $\alpha$ while another one is $\frac{3}{t}$. However, this only  difference leads to great different convergence rates results being  proved. For \eqref{ODE-Nes},  it can be proved that   $F(x(t))-\min_{x} F(x)\sim O(\frac{1}{t^2})$; but for \eqref{ODE-PIAG}, we even have no $O(\frac{1}{t})$ asymptotical rate.

\textbf{The missing constraint:} Only with \eqref{ODE-PIAG}, we can barely obtain any rate. A natural idea to fix this point is adding more information of PIGD to revise the dynamical system \eqref{ODE-PIAG}.
We notice that some important relation between $\ddot{x}(t)$ and $\dot{x}(t)$ is missing.  In the discretization, $\ddot{x}(t)$ is replaced by $\frac{x^{k+1}-2x^{k}+x^{k-1}}{h^2}$, where $h$ is stepsize for discretization. Then it holds that
\begin{equation}
    \left\|\frac{x^{k+1}-2x^{k}+x^{k-1}}{h^2}\right\|\leq \frac{1}{h}\cdot\left(\left\|\frac{x^{k+1}-x^{k}}{h}\right\|+\left\|\frac{x^{k}-x^{k-1}}{h}\right\|\right).
\end{equation}
 Note that both $\frac{x^{k+1}-x^{k}}{h}$ and $\frac{x^{k+1}-x^{k}}{h}$  can be viewed as the discretization of $\dot{x}(t)$. Motivated by this observation, we propose to modify (\ref{ODE-PIAG}) by adding the following constraint
 \begin{equation}\label{constraint}
    \|\ddot{x}(t)\|\leq \theta\|\dot{x}(t)\|,
 \end{equation}
 where $\theta>0$. In Section 2, we study the systems (\ref{ODE1})+(\ref{constraint}) and (\ref{ODE-PIAG})+(\ref{constraint}). With the extra constraint (\ref{constraint}), the sublinear  asymptotical convergence rate can be established. The analysis enables the non-ergodic sublinear convergence rate for heavy-ball  algorithm and PIGD.
\subsection{Related Works}
The inertial term was first proposed in the heavy-ball algorithm \cite{polyak1964some}. When the objective function is twice continuously differentiable, strongly convex (almost quadratic), the Heavy-ball method is proved to  converge linearly. Under weaker assumption that the gradient of the objective function is Lipschitz continuous, \cite{zavriev1993heavy} proved the convergence to a critical point, yet without specifying the convergence rate. The smoothness of objective function is critical for the heavy-ball to converge. In fact, there is an example that the heavy-ball method  diverges for a strongly convex but nonsmooth function \cite{lessard2016analysis}.   Different from the classical gradient methods, heavy-ball algorithm fails to generate a Fej\'{e}r monotone sequence. In general convex and smooth case, the only convergence rate result is ergodic $O(1/k)$ in terms of the function values \cite{ghadimi2015global}. The stochastic inertial parameters are introducing in \cite{loizou2017linearly,loizou2017momentum}.

The iPiano combines heavy-ball method with the proximal mapping as in forward-backward splitting. In the nonconvex case, convergence of the algorithm was thoroughly discussed \cite{ochs2014ipiano}.   The local linear convergence of iPiano and heavy-ball method has been proved in \cite{ochs2016local}. In the strongly convex case, the linear convergence was proved for iPiano  with fixed $\beta_k$ \cite{ochs2015ipiasco}. In the paper \cite{pock2016inertial}, inertial Proximal Alternating Linearized Minimization (iPALM) was introduced as a variant of iPiano for solving two-block  regularized problem. Without the inertial terms, this algorithm reduces to the Proximal Alternating Linearized Minimization (PALM) \cite{bolte2014proximal}, being equivalent to the two-block case of the Coordinate Descent (CD) algorithm \cite{wright2015coordinate}.

\subsection{Contribution and Organization}
In this paper, we present the \emph{first} non-ergodic $O(1/k)$ convergence rate result for PIGD in general convex case. Compared with results in \cite{ochs2015ipiasco}, our convergence is established with a much larger stepsize under the coercive assumption. If the function fails to be coercive, we can choose asymptotic stepsizes. We also present the linear convergence under an error bound condition without assuming strong convexity. Similar to the coercive case, our results hold for relaxed stepsizes. In addition, we extend our result to the coordinate descent version of PIGD. Both cyclic and stochastic index selection strategies are considered. The contributions of this paper are summarized as follows:

1. \textbf{A novel dynamical interpretation:} We propose a novel dynamical system for PIGD, from which we derive the sublinear asymptotical  convergence rate with a proper Lyapunov function under smooth condition.

2. \textbf{The non-ergodic sublinear convergence rate:} We are the first to prove the non-ergodic convergence rates of PIGD. The linear convergence rate is also proved for the objective function without strong convexity. The main idea of proof is to bound the Lyapunov function, and connect this bound to the successive difference of the Lyapunov function.

3. \textbf{Better linear convergence:} Stronger linear convergence results are proved for PIGD. Compared with that in the literature, we proved that PIGD can enjoy relaxed stepsize and inertial parameters. The strong convexity assumption can be weaken. More importantly, we show that the stepsize can be chosen independent of the strong convexity constant.

4. \textbf{Extensions to multi-block version:} The convergence of  multi-block versions of PIGD is studied. Both cyclic and stochastic index selection strategies are considered. The sublinear  and linear convergence rates are proved for  both algorithms.

The rest of the paper is organized as follows. In Section 2, we study the modified dynamical system and present technical lemmas. In Section 3, we show the convergence rates for PIGD. We extend the results to  multi-block deterministic coordinate version of PIGD in Section 4  and to the stochastic coordinate version in Section 5. Section
6 concludes this paper.
\section{Dynamical Motivation and Technical Lemmas}
Due to that the proof is changed greatly without smoothness of the objective function.
Thus, in this part, we first analyze the performance of the modified dynamical system (\ref{ODE1})+(\ref{constraint}), and then consider system (\ref{ODE-PIAG})+(\ref{constraint}).  The existences of the two systems are beyond the scope of this paper and will not be discussed. On another hand, we just want to use the dynamics to  motivate the Lyapunov function for Heavy-ball and PIGD; the existence of the system contributes nothing.  After the analysis of the dynamics,  two necessary lemmas are introduced for the following part of this paper.
\subsection{Performance of   Modified Heavy-ball Dynamical System}
Let us   consider the Lyapunov function as
\begin{align}\label{lyas}
\xi_{f}(t):= f(x(t))+\frac{1}{2}\|\dot{x}(t)\|^2-\min f.
\end{align}
With direct computation, it holds that
\begin{align}\label{system-t1}
\dot{\xi_f}(t)=\langle\nabla f(x(t)),\dot{x}(t)\rangle+\langle\ddot{x}(t),\dot{x}(t)\rangle=-\alpha\|\dot{x}(t)\|^2.
\end{align}
 Assume that $f$ is coercive, noting $\xi_f(t)$ is decreasing and nonnegative, $x(t)$ must be bounded. With the continuity of $\nabla f$, $\nabla f(x(t))$ is also bounded. That means $\ddot{x}(t)+\alpha \dot{x}(t)$ is also bounded. If $\alpha>\theta$, with the triangle inequality,
 \begin{align}
 \|\ddot{x}(t)+\alpha \dot{x}(t)\|\geq \alpha\|\dot{x}(t)\|-\|\ddot{x}(t)\|\geq (\alpha-\theta)\|\dot{x}(t)\|.
 \end{align}
 We then obtain the boundedness of $\ddot{x}(t)$ and $\dot{x}(t)$. Let $x^*\in \textrm{arg}\min f$, we have
 \begin{align}
&f(x(t))-f(x^*)\leq \langle \nabla f(x(t)), x(t)-x^*\rangle\nonumber\\
&\quad\quad\leq \|\nabla f(x(t))\|\cdot\|x(t)-x^*\|\nonumber\\
&\quad\quad\leq (\alpha+\theta)\|\dot{x}(t)\|\cdot\|x(t)-x^*\|.
 \end{align}
 With the boundedness, denote that
 $$R:=\sup_{t\geq 0}\left[\max\{(\alpha+\theta)\cdot\|x(t)-x^*\|,\frac{\|\dot{x}(t)\|}{2}\}\right]<+\infty.$$
 Then we can easily have
 \begin{align}\label{system-t2}
 \xi_{f}(t)^2\leq R^2\|\dot{x}(t)\|^2.
 \end{align}
 With (\ref{system-t1}) and (\ref{system-t2}), we derive
  \begin{align*}
 \xi_{f}(t)^2\leq -\frac{R^2}{\alpha}\dot{\xi_{f}}(t).
 \end{align*}
 That is also
   \begin{align}\label{moto}
 -\frac{\alpha}{R^2}d t\leq\frac{d\xi_f}{\xi_f^2}.
 \end{align}
 Taking integrations of both sides, we then have
 \begin{align*}
f(x(t))-f(x^*)\leq\xi_f(t)\leq \frac{1}{\frac{\alpha}{R^2}t+\xi_f(0)}.
 \end{align*}
 If we just consider (\ref{ODE1}), only convergence can be proved without sublinear  asymptotical rates.  Obviously, (\ref{constraint}) is crucial for the analysis. The analysis above is very easy and only basic calculus is used.
 \subsection{Performance of  PIGD Dynamical System}
For the nonsmooth case, we analyze the following Lyapunov function
\begin{align}\label{lyasn}
\xi_{F}(t):= F(x(t))+\frac{1}{2}\|\dot{x}(t)\|^2-\min F.
\end{align}
Due to that $g$ may be nonsmooth, we cannot directly take differentials of $\xi_{F}(t)$. Instead, we study $\xi_{F}(t+\Delta t)-\xi_{F}(t)$ for small enough $\Delta t>0$. We first review a theoretical result about directional derivative.
\begin{lemma}\cite{rockafellar2015convex}\label{directsub}
For any convex function $f$ and any $x, p \in \mathbb{R}^n$, the directional derivative
$\lim_{t\rightarrow 0+}(f(x + sp) - f(x))/s $ exists, and can be evaluated as
$$\lim_{s\rightarrow 0+}\frac{f(x + sp) - f(x)}{s} = \sup_{\xi\in \partial f(x)}\langle \xi, p \rangle.$$
\end{lemma}
Using Lemma \ref{directsub}, we have the approximation
\begin{align}\label{non-t1}
    &F(X(t+\Delta t))= F(X(t) + \Delta t \dot X(t) + o(\Delta t)) \nonumber\\
    &= F(X(t) + \Delta t \dot X(t)) + o(\Delta t)\nonumber\\
    &=F(X(t))+\langle G(X(t),\dot{X}(t)), \dot{X}(t)\rangle\Delta t+o(\Delta t).
\end{align}
Noting that $\frac{d\|\dot{X}(t)\|^2}{dt}=2\langle \dot{X}(t),\ddot{X}(t)\rangle$, thus
\begin{align}\label{non-t2}
 \frac{1}{2}\|\dot{X}(t+\Delta t)\|^2- \frac{1}{2}\|\dot{X}(t)\|^2=\langle \dot{X}(t),\ddot{X}(t)\rangle\Delta t+o(\Delta t).
\end{align}
With \eqref{non-t1} and \eqref{non-t2}, we are then led to
\begin{align*}
&\xi_{F}(t+\Delta t)-\xi_{F}(t)=\langle G(X(t),\dot{X}(t)), \dot{X}(t)\rangle\Delta t\\
&+\langle \dot{X}(t),\ddot{X}(t)\rangle\Delta t+o(\Delta t)=-\alpha\|\dot{X}(t)\|^2+o(\Delta t)\leq o(\Delta t).
\end{align*}
Thus, we come to the conclusion
$$\limsup_{\Delta t\rightarrow 0+}\frac{\xi_{F}(t+\Delta t)-\xi_{F}(t)}{\Delta t} \le 0.$$
The continuity of $\xi_{F}$ indicates that $\xi_{F}(t)$ is a non-increasing function of $t$.
\subsection{Technical Lemmas}
This parts contains two lemmas on  nonnegative sequences: Lemma \ref{bc} is used to derive the convergence rate. It can be regarded as the discrete form of (\ref{moto}); Lemma \ref{app2} is developed to bound the sequence when inertial parameters are decreasing.
\begin{lemma}[Lemma 3.8, \cite{beck2015convergence}]\label{bc}
Let $\{\alpha_k\}_{k\geq 1}$ be nonnegative sequence of real numbers satisfying
\begin{equation*}
    \alpha_k-\alpha_{k+1}\geq \gamma\alpha_{k+1}^2.
\end{equation*}
Then we have
\begin{equation*}
    \alpha_k=O(\frac{1}{k}).
\end{equation*}
\end{lemma}

\begin{lemma}\label{app2}
Let $\{t_k\}_{k\geq 0}$ be a nonnegative sequence and follow the condition
\begin{equation}\label{app2-assump}
    t_{k+1}\leq (1+\beta_k)t_k+\beta_k t_{k-1}.
\end{equation}
If $\{\beta\}_{k\geq 0}$ is descending and
$$\sum_{k}\beta_k<+\infty,$$
$\{t_k\}_{k\geq 0}$ is bounded.
\end{lemma}

\section{Convergence Rates}
In this section, we prove  convergence rates of PIGD. The core of the proof is to construct a proper Lyapunov function. In fact, there are various Lyapunov functions for PIGD.
We first present an obvious one in the following lemma.

\begin{lemma}\label{le1}
Suppose $f$ is convex and has $L$-Lipschitz gradient, and $g$ is convex, and  $\min F>-\infty$. Let $\{x^k\}_{k\geq 0}$ be generated by PIGD  with non-increasing $\{\beta_k\}_{k\geq 0}\subseteq [0,1)$. Choosing the step size
 $$\gamma_k=\frac{2(1-\beta_k)c}{L}$$
 for arbitrary fixed $0<c<1$, we have
\begin{align}\label{nresult}
    & \left[F(x^{k})+\frac{\beta_k}{2\gamma_k}\|x^{k}-x^{k-1}\|^2\right]\nonumber\\
    &\quad\quad-\left[F(x^{k+1})+\frac{\beta_{k+1}}{2\gamma_{k+1}}\|x^{k+1}-x^{k}\|^2\right]\nonumber\\
     &\geq(\frac{1-\beta_k}{\gamma_k}-\frac{L}{2})\|x^{k+1}-x^k\|^2.
\end{align}
\end{lemma}

However, the  Lyapunov function in Lemma \ref{le1} cannot derive the non-ergodic sublinear rate. This is because we cannot derive an inequality like \eqref{system-t2} for this Lyapunov function.
Thus  we employ the a modified  Lyapunov function which reads as
\begin{align}\label{Lyapunov}
   \xi_k:= F(x^k)+\delta_k\|x^k-x^{k-1}\|^2-\min F,
\end{align}
where
\begin{equation}\label{eq:delta}
    \delta_k:=\frac{1}{2}(\frac{1}{\gamma_k}-\frac{L}{2}).
\end{equation}
%
We present a very useful technique lemma which is the key to results.
\begin{lemma}\label{lem-sub}
Suppose the conditions of Lemma \ref{le1} hold.
 Let $\overline{x^k}$ denote the projection of $x^k$ onto $\emph{arg}\min{F}$, assumed to exist, and define
\begin{align}\label{eq:alpha}
 \varepsilon_k: =\frac{4c\delta_{k+1}^2}{(1-c)L}+\frac{4c}{(1-c)L\gamma^2_k}.
\end{align}
Then it holds
\begin{align}\label{result-sub}
   (\xi_{k+1})^2&\leq\varepsilon_k\times(\xi_k-\xi_{k+1})\nonumber\\
   &\times(2\|x^{k+1}-\overline{x^{k+1}}\|^2+\|x^{k+1}-x^{k}\|^2).
\end{align}
\end{lemma}

The inequality \eqref{result-sub}  can be regarded as the discretization of  \eqref{system-t2}.
\subsection{Sublinear Convergence Rate under General Convexity}

In this subsection, we present the sublinear of the convex PIGD. The coercivity of the function is critical for the analysis.
If $F$ is coercive, the parameter $\beta_k$ can be bounded from $0$; however, if $F$ fails to be promised to be coercive, $\beta_k$ must be descending to zero. Thus, this subsection will be divided into two parts in  term of the coercivity.

\subsubsection{$F$ is Coercive}
First, we present the non-ergodic $O(\frac{1}{k})$ convergence rate of the function value. The rate can be derived if $\{\beta_k\}_{k\geq 0}$ is bounded from $0$ and $1$.
\begin{theorem}\label{sub-n}
Assume the conditions of Lemma \ref{le1} hold, and
 $$0<\inf_k\beta_k\leq\beta_k\leq\beta_0<1.$$
Then we have
\begin{equation}\label{sub-n-result}
    F(x^k)-\min F= O(\frac{1}{k}).
\end{equation}
\end{theorem}

To the best of our knowledge, this is the first time to prove the non-ergodic $O(\frac{1}{k})$ convergence rate in the perspective of function values for PIGD and heavy-ball method in the convex case.
\subsubsection{$F$ Fails to Be Coercive}
In this case, to obtain the boundedness of the sequence $\{x^k\}_{k\geq 0}$, we must employ diminishing $\beta_k$, i.e., $\lim_k\beta_k=0$. The following lemma can derive the needed boundedness.

\begin{lemma}\label{lem-b}
Suppose the conditions of Lemma \ref{le1} hold and
 $$\beta_k=\frac{1}{(k+1)^{\theta}},$$
 where $\theta>1$.
Let $\{x^k\}_{k\geq 0}$ be generated by PIGD, then $\{x^k\}_{k\geq 0}$ is bounded.
\end{lemma}

Now, we are prepared to present the $O(1/k)$ rate of the function values when $F$ is not coercive.
\begin{theorem}\label{th-b}
Suppose the conditions of Lemma \ref{lem-b} hold.
Let $\{x^k\}_{k\geq 0}$ be generated by PIGD, then  we have
\begin{equation*}\label{th-b-result}
    F(x^k)-\min F=O(\frac{1}{k}).
\end{equation*}
\end{theorem}

\subsection{Linear Convergence   under Optimal Strong
Convexity Condition}
We say that the function $F$ satisfies the optimal strong
convexity condition if
\begin{equation}\label{condtion}
    F(x)-\min F\geq \nu\|x-\overline{x}\|^2,
\end{equation}
where $\overline{x}$ is the projection of $x$ onto the set $\textrm{arg}\min F$, and $\nu>0$. This condition is much weaker than the strong convexity.
\begin{theorem}\label{thm-linear}
Suppose the conditions of Theorem \ref{sub-n} hold, and $F$ satisfies (\ref{condtion}).
 Then we have
\begin{equation*}\label{thm-linear-result}
   F(x^k)-\min F= O(\omega^k),
\end{equation*}
for some $0<\omega<1$.
\end{theorem}
Compared with previous linear convergence result presented in \cite{ochs2015ipiasco}. Our result enjoys three advantages: 1. The strongly convex assumption is weaken to (\ref{condtion}). 2. More general parameters setting can be used. 3. The stepsizes and inertial parameters are independent with the strongly convex constants.

\section{Cyclic Coordinate PIGD}
This part analyzes the cyclic coordinate inertial proximal algorithm. The two-block version is proposed in \cite{pock2016inertial}, which focuses on the nonconvex case. Here, we consider the multi-block version  and prove its convergence rate under convexity assumption. The minimization problem can be described as
\begin{equation}\label{multi}
    \min_{x_1,x_2,\ldots,x_m}\left\{\underbrace{f(x_1,x_2,\ldots,x_m)+\sum_{i=1}^m g_i(x_i)}_{:=F(x_1,x_2,\ldots,x_m)}\right\},
\end{equation}
where $f$ and $g_i$ ($i=1,\ldots,m$) are all convex. We use the notation
$$\nabla_i^k f:=\nabla_i f(x_1^{k+1},\ldots,x_{i-1}^{k+1},x_{i}^{k}\ldots,x_m^k),~ x^k:=(x^k_1,x^k_2,\ldots,x^k_m).$$
 The cyclic coordinate descent inertial algorithm runs as: for $i$ from $1$ to $m$,
\begin{equation}\label{schemec}
    x^{k+1}_i=\textbf{prox}_{\gamma_{k,i} g_i}[x^k_i-\gamma_{k,i}\nabla_i^k f+\beta_{k,i}(x^k_i-x^{k-1}_{i})],
\end{equation}
where $\gamma_{k,i},\beta_{k,i}>0$. The iPALM can be regarded as the two-block case of this algorithm. The function $f$ is assumed to satisfy
\begin{align}\label{picl}
   &\|\nabla_i f(x_1,x_2,\ldots,x_i^1,\ldots,x_m)\nonumber\\
   &-\nabla_i f(x_1,x_2,\ldots,x_i^2,\ldots,x_m)\|\leq L_i\|x_i^1-x_i^2\|
\end{align}
for any $x_1,x_2,\ldots,x_{i-1},x_{i+1},\ldots,x_m$, and $x_i^1, x_i^2$, and $i\in[1,2,\ldots,m]$. With (\ref{picl}), we can easily obtain
 \begin{align}\label{lpm}
   & f(x_1,x_2,\ldots,x_i^1,\ldots,x_m)\leq   f(x_1,x_2,\ldots,x_i^2,\ldots,x_m)\nonumber\\
    &+\langle \nabla_i f(x_1,x_2,\ldots,x_i^2,\ldots,x_m),x_i^1-x_i^2\rangle+\frac{L_i}{2}\|x_i^1-x_i^2\|^2.
 \end{align}
 The proof is similar to [Lemma 1.2.3,\cite{nesterov2013introductory}] and will not be reproduced.
 In the following part of this paper, we use the following assumption

\textbf{A1}: for any $i\in[1,2,\ldots,m]$, the sequence $(\beta_{k,i})_{k\geq 0}\subseteq [0,1)$ is non-increasing.

\begin{lemma}\label{md-le1}
 Let $f$ be a convex function  satisfying (\ref{picl})  and $g_i$ is convex ($i\in[1,2,\ldots,m]$), and finite $\min F$. Assume  $\{x^k\}_{k\geq 0}$ is generated by scheme (\ref{schemec}) and assumption \textbf{A1} is satisfied. Choose the step size
 $$\gamma_{k,i}=\frac{2(1-\beta_{k,i})c}{L_i},~\quad i\in[1,2,\ldots,m]$$
 for arbitrary fixed $0<c<1$. Then we can obtain
\begin{align}\label{md-result}
    & \left[F(x^{k})+\sum_{i=1}^m\frac{\beta_{k,i}}{2\gamma_{k,i}}\|x^{k}_i-x^{k-1}_i\|^2\right]\nonumber\\
    &\quad\quad-\left[F(x^{k+1})+\sum_{i=1}^m\frac{\beta_{k+1,i}}{2\gamma_{k+1,i}}\|x^{k+1}_i-x^{k}_i\|^2\right]\nonumber\\
     &\geq\frac{(1-c)\underline{L}}{2c}\|x^{k+1}-x^k\|^2,
\end{align}
where $\underline{L}=\min_{i\in[1,2,\ldots,m]}\{L_i\}$.
\end{lemma}

The following Lyapunov function is used for cyclic coordinate PIGD
\begin{align}\label{Lyapunov-m}
   \hat{\xi}_k:= F(x^k)+\sum_{i=1}^m \delta_{k,i}\|x^k_i-x^{k-1}_i\|^2-\min F,
\end{align}
where
\begin{equation}\label{eq:delta-m}
    \delta_{k,i}:=\frac{1}{2}(\frac{1}{\gamma_{k,i}}-\frac{L_i}{2}).
\end{equation}
With this  Lyapunov function, we can present the following lemma.
\begin{lemma}\label{lem-sub-m}
Suppose the conditions Lemma \ref{md-le1} hold.
 Let $\overline{x^k}$ denote the projection of $x^k$ onto $\emph{arg}\min{F}$, assumed to exist, and define
\begin{align*}
 \hat{\varepsilon}_k: =\max\{\frac{4c}{(1-c)\underline{L}}\sum_{i=1}^m\left(\delta_{k+1,i}^2+L_i^2\right),\frac{4c}{(1-c)\underline{L}}\sum_{i=1}^m\frac{1}{\gamma^2_{k,i}}\}.
\end{align*}
Then it holds  that
\begin{align}\label{result-sub-m}
   (\hat{\xi}_{k+1})^2&\leq\hat{\varepsilon}_k\times(\hat{\xi}_k-\hat{\xi}_{k+1})\nonumber\\
   &\times(3\|x^{k+1}-\overline{x^{k+1}}\|^2+\|x^{k-1}-x^{k}\|^2).
\end{align}
\end{lemma}

\subsection{Sublinear Convergence Rates of Cyclic Coordinate PIGD}
This part proves the sublinear convergence rates of  cyclic coordinate descent inertial algorithm. In multi-block case, it is always to assume that the objective function is coercive.  Like the previous section, we  obtain the $O(1/k)$ convergence rate of the algorithm if $F$ is coercive.

\begin{theorem}\label{them-sub-m}
Suppose the conditions of Lemma \ref{md-le1} hold, $F$ is coercive and
 $$0<\inf_k\beta_k\leq\beta_{k,i}\leq\beta_0<1,~i\in[1,2,\ldots,m].$$
Then we have
\begin{align}
   F(x^k)-\min F= O(\frac{1}{k}).
\end{align}
\end{theorem}
\subsection{Linear Convergence Rate of Cyclic Coordinate   PIGD}
In this part, we establish the linear convergence result of cyclic coordinate   PIGD.
\begin{theorem}\label{them-linear-m}
Suppose the conditions of Lemma \ref{md-le1} hold, $F$ satisfies (\ref{condtion}), and
 $$0<\inf_k\beta_k\leq\beta_{k,i}\leq\beta_0<1,~i\in[1,2,\ldots,m].$$
Then we have
\begin{align*}
   F(x^k)-\min F= O(\omega^k)
\end{align*}
for some $0<\omega<1$.
\end{theorem}

\section{Stochastic Coordinate PIGD}
We still aim to minimizing problem (\ref{multi}) but using the stochastic index selection strategy. In the $k$-th iteration, pick $i_k$ uniformly from $[1,2,\ldots,m]$, and then update
\begin{equation}\label{schemes}
    x^{k+1}_{i_k}=\textbf{prox}_{\gamma_{k} g_{i_k}}[x^k_{i_k}-\gamma_{k}\nabla_{i_k} f(x^k)+\beta_{k}(x^k_{i_k}-x^{k-1}_{i_k})].
\end{equation}
The sub-algebra $\chi^k$ is defined as
\begin{equation}\label{sub-al}
    \chi^k:=\sigma(x^0,x^1,\ldots,x^k).
\end{equation}
In this section, we use the following assumption
following assumption

\textbf{A2}:  the sequence $(\beta_{k})_{k\geq 0}\subseteq [0,1)$ is non-increasing.

Denote that
\begin{equation}\label{DS}
    S_{\gamma}(x)=x-\textbf{prox}_{\gamma g}[x-\gamma\nabla f(x)].
\end{equation}
It is easy to see that
\begin{equation}\label{equal}
    S_{\gamma}(x^*)=\textbf{0}\Longleftrightarrow x^*~~\textrm{minimize}~~F(x)
\end{equation}
for any $\gamma>0$. We first present a general convergence result of  stochastic coordinate PIGD. The $O(1/k)$ convergence rate of the algorithm  is proved for the successive difference of the points.
\begin{lemma}\label{ms-le1}
 Let $f$ be a convex function  whose gradient is Lipschitz continuous with $L$, and $g_i$ is convex ($i\in[1,2,\ldots,m]$), and finite $\min F$. Assume $\{x^k\}_{k\geq 0}$ is generated by scheme (\ref{schemes}) and assumption \textbf{A2} is satisfied. Choose the step size
 $$\gamma_{k}=\frac{2(1-\beta_{k}/\sqrt{m})c}{L}$$
 for arbitrary fixed $0<c<1$. Then  we can obtain
\begin{align}\label{ms-result}
    & \left[\EE F(x^{k})+\frac{\beta_k}{2\sqrt{m}\gamma_k}\EE\|x^{k}-x^{k-1}\|^2\right]\nonumber\\
    &\quad\quad-\left[\EE F(x^{k+1})+\frac{\beta_{k+1}}{2\sqrt{m}\gamma_{k+1}}\EE\|x^{k+1}-x^{k}\|^2\right]\nonumber\\
     &\geq(\frac{1-\beta_k/\sqrt{m}}{\gamma_k}-\frac{L}{2})\EE\|x^{k+1}-x^k\|^2.
\end{align}
And we further have
\begin{equation}\label{ms-result-2}
   \EE( \min_{0\leq i\leq k}\|x^{i+1}-x^i\|^2)=o(\frac{1}{k}).
\end{equation}
\end{lemma}

\begin{theorem}\label{ms-th1}
Suppose that the conditions of Lemma \ref{ms-le1} hold  and $0\leq \beta<\sqrt{m}$.
Then we have
\begin{equation}\label{ms-th-result}
    \EE(\min_{0\leq i\leq k}\|S_{\gamma}(x^k)\|^2)=o(\frac{1}{k}).
\end{equation}
\end{theorem}

The non-ergodic sublinear rate of stochastic coordinate PIGD is hard to prove with previous proof techniques and routines. This is because the bound of $\EE F(x^{k+1})-\min F$ can barely be estimated for the conditional expectation rules\footnote{ If $g\equiv 0$, the non-ergodic sublinear convergence rate can be proved. Details can be founded in \cite{sun2018non}.}. Therefore, the non-ergodic convergence of stochastic coordinate PIGD will not be presented.  But the linear convergence rate can be derived provided we choose proper stepsizes and inertial parameter.
We now introduce the linear convergence results of  stochastic coordinate PIGD.
\begin{lemma}\label{lem:linear-convergence}
Assume the function $F$ satisfies the optimal strong convexity condition (\ref{condtion}), and $\{x^k\}_{k\geq 0}$ is generated by the scheme (\ref{schemes}).  By choosing
$$\beta_k\equiv\beta=\frac{\gamma\nu}{4m},$$
 we have
\begin{align*}
&\underline{\ell}\EE[\|x^{k}-\overline{x^{k}}\|^2+2\gamma(F(x^{k})-F^*)] \\
&\quad\leq\EE[\|x^{k}-\overline{x^{k}}\|^2+2\gamma(F(x^{k})-F^*)]\nonumber\\
&-\EE[\|x^{k+1}-\overline{x^{k+1}}\|^2+2\gamma(F(x^{k+1})-F^*)]\nonumber\\
&\quad-(1-\gamma L-2\beta)\EE\|x^{k+1}-x^k\|^2+\frac{\beta}{m}\EE\|x^{k}-x^{k-1}\|^2,\nonumber
\end{align*}
where $\underline{\ell}:=\frac{\min\{\nu,1\}}{2m}\gamma$ and $\overline{x^k}$ is the projection of $x^k$ onto $\emph{arg}\min F$.
\end{lemma}

To present the linear convergence rate, we introduce a constant. We consider the following equation
$$\frac{\min\{\nu,1\}\nu}{8m^3}\gamma^2+\left(L+\frac{\nu}{2m}-\frac{\nu}{4m^2}\right)\gamma-1=0.$$
It is easy to see that there exist a positive root in $(0,\frac{1}{L})$, which is denoted as $\gamma_0$ here.
\begin{theorem}\label{th:linear-convergence}
Suppose the conditions of Lemma \ref{lem:linear-convergence} hold. If
$$0<\gamma<\gamma_0,~\,~\beta=\frac{\gamma\nu}{4m},$$
 we have
\begin{align*}
\EE[F(x^{k})-F^*] =O\left((1-\frac{\gamma\nu}{2m})^k\right).
\end{align*}
\end{theorem}

Noticing  that $0<\gamma<\gamma_0<\frac{1}{L}$ and $\nu\leq L$,  $0<\frac{\gamma\nu}{2m}<\frac{1}{2m}$.
\section{Conclusion}
In this paper, we focus on the non-ergodic complexity of the  proximal inertial gradient descent in the convex setting. We prove the non-ergodic  sublinear convergence  rate of the algorithm and linear convergence rate under larger stepsize. We extend our results to the multi-block inertial algorithm. For both cyclic and stochastic index selection strategies, the convergence rates  are proved.

%
%


\begin{thebibliography}{10}

\bibitem{ochs2014ipiano}
P.~Ochs, Y.~Chen, T.~Brox, and T.~Pock, ``ipiano: Inertial proximal algorithm
  for nonconvex optimization,'' {\em SIAM Journal on Imaging Sciences}, vol.~7,
  no.~2, pp.~1388--1419, 2014.

\bibitem{combettes2005signal}
P.~L. Combettes and V.~R. Wajs, ``Signal recovery by proximal forward-backward
  splitting,'' {\em Multiscale Modeling \& Simulation}, vol.~4, no.~4,
  pp.~1168--1200, 2005.

\bibitem{polyak1964some}
B.~T. Polyak, ``Some methods of speeding up the convergence of iteration
  methods,'' {\em USSR Computational Mathematics and Mathematical Physics},
  vol.~4, no.~5, pp.~1--17, 1964.

\bibitem{su2014differential}
W.~Su, S.~Boyd, and E.~Candes, ``A differential equation for modeling
  nesterov's accelerated gradient method: Theory and insights,'' in {\em
  Advances in Neural Information Processing Systems; Journal of Machine
  Learning Research}, pp.~2510--2518, 2014.

\bibitem{sun2017asynchronous}
T.~Sun, R.~Hannah, and W.~Yin, ``Asynchronous coordinate descent under more
  realistic assumptions,'' {\em NIPS}, 2017.

\bibitem{wibisono2016variational}
A.~Wibisono, A.~C. Wilson, and M.~I. Jordan, ``A variational perspective on
  accelerated methods in optimization,'' {\em proceedings of the National
  Academy of Sciences}, vol.~113, no.~47, pp.~E7351--E7358, 2016.

\bibitem{attouch2018fast}
H.~Attouch, Z.~Chbani, J.~Peypouquet, and P.~Redont, ``Fast convergence of
  inertial dynamics and algorithms with asymptotic vanishing viscosity,'' {\em
  Mathematical Programming}, vol.~168, no.~1-2, pp.~123--175, 2018.

\bibitem{alvarez2000minimizing}
F.~Alvarez, ``On the minimizing property of a second order dissipative system
  in hilbert spaces,'' {\em SIAM Journal on Control and Optimization}, vol.~38,
  no.~4, pp.~1102--1119, 2000.

\bibitem{rockafellar2015convex}
R.~T. Rockafellar, {\em Convex analysis}.
\newblock Princeton university press, 2015.

\bibitem{beck2009fast}
A.~Beck and M.~Teboulle, ``A fast iterative shrinkage-thresholding algorithm
  for linear inverse problems,'' {\em SIAM journal on imaging sciences},
  vol.~2, no.~1, pp.~183--202, 2009.

\bibitem{nesterov2013introductory}
Y.~Nesterov, {\em Introductory lectures on convex optimization: A basic
  course}, vol.~87.
\newblock Springer Science \& Business Media, 2013.

\bibitem{zavriev1993heavy}
S.~Zavriev and F.~Kostyuk, ``Heavy-ball method in nonconvex optimization
  problems,'' {\em Computational Mathematics and Modeling}, vol.~4, no.~4,
  pp.~336--341, 1993.

\bibitem{lessard2016analysis}
L.~Lessard, B.~Recht, and A.~Packard, ``Analysis and design of optimization
  algorithms via integral quadratic constraints,'' {\em SIAM Journal on
  Optimization}, vol.~26, no.~1, pp.~57--95, 2016.

\bibitem{ghadimi2015global}
E.~Ghadimi, H.~R. Feyzmahdavian, and M.~Johansson, ``Global convergence of the
  heavy-ball method for convex optimization,'' in {\em Control Conference
  (ECC), 2015 European}, pp.~310--315, IEEE, 2015.

\bibitem{loizou2017linearly}
N.~Loizou and P.~Richt{\'a}rik, ``Linearly convergent stochastic heavy ball
  method for minimizing generalization error,'' {\em arXiv preprint
  arXiv:1710.10737}, 2017.

\bibitem{loizou2017momentum}
N.~Loizou and P.~Richt{\'a}rik, ``Momentum and stochastic momentum for
  stochastic gradient, newton, proximal point and subspace descent methods,''
  {\em arXiv preprint arXiv:1712.09677}, 2017.

\bibitem{ochs2016local}
P.~Ochs, ``Local convergence of the heavy-ball method and ipiano for non-convex
  optimization,'' {\em arXiv preprint arXiv:1606.09070}, 2016.

\bibitem{ochs2015ipiasco}
P.~Ochs, T.~Brox, and T.~Pock, ``ipiasco: Inertial proximal algorithm for
  strongly convex optimization,'' {\em Journal of Mathematical Imaging and
  Vision}, vol.~53, no.~2, pp.~171--181, 2015.

\bibitem{pock2016inertial}
T.~Pock and S.~Sabach, ``Inertial proximal alternating linearized minimization
  (ipalm) for nonconvex and nonsmooth problems,'' {\em SIAM Journal on Imaging
  Sciences}, vol.~9, no.~4, pp.~1756--1787, 2016.

\bibitem{bolte2014proximal}
J.~Bolte, S.~Sabach, and M.~Teboulle, ``Proximal alternating linearized
  minimization for nonconvex and nonsmooth problems,'' {\em Mathematical
  Programming}, vol.~146, no.~1-2, pp.~459--494, 2014.

\bibitem{wright2015coordinate}
S.~J. Wright, ``Coordinate descent algorithms,'' {\em Mathematical
  Programming}, vol.~151, no.~1, pp.~3--34, 2015.

\bibitem{beck2015convergence}
A.~Beck, ``On the convergence of alternating minimization for convex
  programming with applications to iteratively reweighted least squares and
  decomposition schemes,'' {\em SIAM Journal on Optimization}, vol.~25, no.~1,
  pp.~185--209, 2015.

\bibitem{sun2018non}
T.~Sun, P.~Yin, D.~Li, C.~Huang, L.~Guan, and H.~Jiang, ``Non-ergodic
  convergence analysis of heavy-ball algorithms,'' {\em AAAI 2019}, 2019.

\bibitem{davis2016convergence}
D.~Davis and W.~Yin, ``Convergence rate analysis of several splitting
  schemes,'' in {\em Splitting methods in communication, imaging, science, and
  engineering}, pp.~115--163, Springer, 2016.

\end{thebibliography}
\appendices
\section{Proof of Lemma \ref{app2}}
Adding $\beta_k t_k$ to both sides of (\ref{app2-assump}),
 \begin{align}\label{app2-t1}
    t_{k+1}+\beta_k t_k&\leq (1+\beta_k)t_k+\beta_k t_{k-1}+\beta_k t_k\nonumber\\
    &\leq (1+2\beta_k)(t_k+\beta_k t_{k-1}).
 \end{align}
 Noting the decent of $\{\beta_k\}_{k\geq 0}$, (\ref{app2-t1}) is actually
  \begin{equation*}
    t_{k+1}+\beta_k t_k\leq (1+2\beta_k)(t_k+\beta_{k-1} t_{k-1}).
 \end{equation*}
 Letting
 $$h_k:=t_k+\beta_{k-1} t_{k-1},$$
 we then have
   \begin{equation*}
    h_{k+1}\leq (1+2\beta_k)h_k\leq e^{2\beta_k}h_k.
 \end{equation*}
 Thus, for any $k$
 \begin{equation*}
    h_{k+1}\leq  e^{2\sum_{i=1}^k\beta_i}h_1<+\infty.
 \end{equation*}
 The boundedness of $\{h_k\}_{k\geq 0}$ directly yields the boundedness of $\{t_k\}_{k\geq 0}$.
 \section{Proof of Lemma \ref{le1}}
The K.K.T. condition of updating $x^{k+1}$ directly gives
\begin{equation}\label{le1-t1}
    \frac{x^k-x^{k+1}}{\gamma_k}-\nabla f(x^k)+\frac{\beta_k}{\gamma_k}(x^k-x^{k-1})\in \partial g(x^{k+1}).
\end{equation}
With the convexity of $g$, we have
\begin{align}\label{le1-t2}
    &g(x^{k+1})-g(x^k)\nonumber\\
    &\quad\leq \langle \frac{x^{k+1}-x^k}{\gamma_k}+\nabla f(x^k)+\frac{\beta_k}{\gamma_k}(x^{k-1}-x^k), x^k-x^{k+1}\rangle.
\end{align}
The Lipschitz continuity of $\nabla f$ tells us
\begin{equation}\label{le1-t3}
    f(x^{k+1})-f(x^k)\leq \langle -\nabla f(x^k),x^k-x^{k+1}\rangle+\frac{L}{2}\|x^{k+1}-x^k\|^2.
\end{equation}
Combining (\ref{le1-t2}) and (\ref{le1-t3}),
\begin{align}\label{le1-t4}
    &F(x^{k+1})-F(x^{k})\nonumber\\
    &\overset{a)}{\leq} \frac{\beta_k}{\gamma_k}\langle x^{k}-x^{k-1}, x^{k+1}-x^k\rangle+(\frac{L}{2}-\frac{1}{\gamma_k})\|x^{k+1}-x^k\|^2\nonumber\\
    &\overset{b)}{\leq}\frac{\beta_k}{2\gamma_k}\|x^{k}-x^{k-1}\|^2+(\frac{L}{2}-\frac{1}{\gamma_k}+\frac{\beta_k}{2\gamma_k})\|x^{k+1}-x^k\|^2.
 \end{align}
where $a)$ comes from (\ref{le1-t2})+(\ref{le1-t3}), and $b)$ uses the Schwarz inequality $\langle x^{k}-x^{k-1}, x^{k+1}-x^k\rangle\leq \frac{1}{2}\|x^{k}-x^{k-1}\|^2+\frac{1}{2}\|x^{k+1}-x^k\|^2$. With direct calculations, we then obtain
\begin{align}\label{le1-t5}
    & \left[F(x^{k})+\frac{\beta_k}{2\gamma_k}\|x^{k}-x^{k-1}\|^2\right]\nonumber\\
    &\quad\quad-\left[F(x^{k+1})+\frac{\beta_{k}}{2\gamma_{k}}\|x^{k+1}-x^{k}\|^2\right]\nonumber\\
     &\quad\quad\geq(\frac{1-\beta_k}{\gamma_k}-\frac{L}{2})\|x^{k+1}-x^k\|^2.
\end{align}
With the non-increasity of $\{\beta_k\}_{k\geq 0}$, $\{\frac{\beta_k}{2\gamma_k}=\frac{\beta_k L}{4(1-\beta_k)c}\}_{k\geq 0}$ is also non-increasing.
Thus,
(\ref{nresult}) can be obtained by \eqref{le1-t5} and
\begin{align*}
    & \left[F(x^{k})+\frac{\beta_k}{2\gamma_k}\|x^{k}-x^{k-1}\|^2\right]\\
    &\quad\quad-\left[F(x^{k+1})+\frac{\beta_{k+1}}{2\gamma_{k+1}}\|x^{k+1}-x^{k}\|^2\right]\\
     &\geq \left[F(x^{k})+\frac{\beta_k}{2\gamma_k}\|x^{k}-x^{k-1}\|^2\right]\\
     &\quad\quad-\left[F(x^{k+1})+\frac{\beta_{k}}{2\gamma_{k}}\|x^{k+1}-x^{k}\|^2\right].
\end{align*}

\section{Proof of Lemma \ref{lem-sub}}
With direct computation and Lemma \ref{le1}, we have
\begin{align}\label{sketch-1}
    &\xi_k- \xi_{k+1}\geq\frac{1}{2}(\frac{1-\beta_k}{\gamma_k}-\frac{L}{2})\nonumber\\
    &\quad\times(\|x^{k+1}-x^k\|^2+\|x^k-x^{k-1}\|^2)\nonumber\\
    &\quad=\frac{L}{4}(\frac{1}{c}-1)\times(\|x^{k+1}-x^k\|^2+\|x^k-x^{k-1}\|^2).
\end{align}
The convexity of $g$ yields
\begin{align*}
    g(x^{k+1})- g(\overline{x^{k+1}})\leq\langle\widetilde{\nabla} g(x^{k+1}),x^{k+1}-\overline{x^{k+1}}\rangle,
\end{align*}
where $\widetilde{\nabla} g(x^{k+1})\in \partial g(x^{k+1})$. With (\ref{le1-t1}), we then have
\begin{align}\label{lem-sub-t4}
    &g(x^{k+1})- g(\overline{x^{k+1}})\nonumber\\
    &\leq\langle     \frac{x^k-x^{k+1}}{\gamma_k}-\nabla f(x^k)+\frac{\beta_k}{\gamma_k}(x^k-x^{k-1}), x^{k+1}-\overline{x^{k+1}}\rangle.
\end{align}
Similarly, we have
\begin{equation}\label{lem-sub-t5}
    f(x^{k+1})-f(\overline{x^{k+1}})\leq \langle\nabla f(x^{k+1}),x^{k+1}-\overline{x^{k+1}}\rangle.
\end{equation}
Summing (\ref{lem-sub-t4}) and (\ref{lem-sub-t5})  yields
\begin{align}\label{lem-sub-t6}
    &F(x^{k+1})-F(\overline{x^{k+1}})\nonumber\\
    &\quad\leq \frac{\beta_k}{\gamma_k}\langle x^{k}-x^{k-1},x^{k+1}-\overline{x^{k+1}}\rangle\nonumber\\
    &\quad+\langle  \frac{x^k-x^{k+1}}{\gamma_k},x^{k+1}-\overline{x^{k+1}}\rangle\nonumber\\
    &\quad\overset{a)}{\leq}\frac{\beta_k}{\gamma_k}\|x^{k}-x^{k-1}\|\cdot\|x^{k+1}-\overline{x^{k+1}}\|\nonumber\\
    &\quad+\frac{1}{\gamma_k}\|  x^k-x^{k+1}\|\cdot\|x^{k+1}-\overline{x^{k+1}}\|\nonumber\\
    &\quad\overset{b)}{\leq}\frac{1}{\gamma_k}\left(\|x^k-x^{k+1}\|+\|x^{k}-x^{k-1}\|\right)\times\|x^{k+1}-\overline{x^{k+1}}\|,
\end{align}
where $a)$ is due to the Schwarz inequalities, $b)$ depends on the fact $0\leq\beta_k<1$. With (\ref{Lyapunov}) and (\ref{lem-sub-t6}), we have
\begin{align}\label{lem-sub-t7}
    \xi_{k+1}&\leq\frac{1}{\gamma_k}\left(\|x^k-x^{k+1}\|+\|x^{k}-x^{k-1}\|\right)\nonumber\\
    &\times\|x^{k+1}-\overline{x^{k+1}}\|+\delta_{k+1}\|x^{k+1}-x^{k}\|^2.\nonumber
\end{align}
Let
 \begin{align}
    a^k:=\left(
                              \begin{array}{c}
                              \frac{1}{\gamma_k}\|x^k-x^{k+1}\| \\
                              \frac{1}{\gamma_k}\|x^{k}-x^{k-1}\| \\
                               \delta_{k+1}\|x^{k+1}-x^{k}\| \\
                              \end{array}
                            \right),\quad b^k:=\left(
                              \begin{array}{c}
                              \|x^{k+1}-\overline{x^{k+1}}\| \\
                              \|x^{k+1}-\overline{x^{k+1}}\| \\
                               \|x^{k+1}-x^{k}\| \\
                              \end{array}
                            \right).
\end{align}
Using this and the definition of $\xi_{k+1}$ (\ref{Lyapunov}), we have
\begin{align}
    (\xi_{k+1})^2=\left|\langle a^k,b^k\rangle\right|^2\leq\|a^k\|^2\times\|b^k\|^2.
\end{align}
Direct calculation yields
\begin{equation*}
    \|a^k\|^2\leq(\delta_{k+1}^2+\frac{1}{\gamma^2_k})\times(\|x^k-x^{k+1}\|^2+\|x^{k-1}-x^{k}\|^2)
\end{equation*}
and
\begin{equation*}
    \|b^k\|^2\leq 2\|x^{k+1}-\overline{x^{k+1}}\|^2+\|x^{k+1}-x^{k}\|^2.
\end{equation*}
Thus we derive
\begin{align}\label{sketch-2}
    (\xi_{k+1})^2&\leq(\delta_{k+1}^2+\frac{1}{\gamma^2_k})\times(\|x^k-x^{k+1}\|^2+\|x^{k-1}-x^{k}\|^2)\nonumber\\
    &\quad\times(2\|x^{k+1}-\overline{x^{k+1}}\|^2+\|x^{k+1}-x^{k}\|^2).
\end{align}
Combining (\ref{sketch-1}) and (\ref{sketch-2}), we then prove the result.

\subsection{Proof of Theorem \ref{sub-n}}
 By Lemma \ref{lem-sub},  $\sup_{k}\{\xi_{k}\}<+\infty$, thus, $\sup_{k}\{F(x^k)\}<+\infty$ and $\sup_{k}\{\|x^{k}-x^{k-1}\|^2\}<+\infty$. Noting the coercivity of $F$, sequences $\{x^k\}_{k\geq 0}$ and $\{\overline{x^k}\}_{k\geq 0}$ are  bounded. With the assumptions on $\gamma_k$ and $\beta_k$, $\sup_k\{\varepsilon_k\}<+\infty$. Thus, $\left\{\varepsilon_k(2\|x^{k+1}-\overline{x^{k+1}}\|^2+\|x^{k-1}-x^{k}\|^2)\right\}_{k\geq 1}$ is bounded; and we assume the bound is $R$, i.e.,
\begin{equation*}
    \sup_{k}\{\varepsilon_k(2\|x^{k+1}-\overline{x^{k+1}}\|^2+\|x^{k-1}-x^{k}\|^2)\}\leq R.
\end{equation*}
By Lemma \ref{lem-sub},  we then have
\begin{equation*}
    \xi_{k+1}^2 \leq R(\xi_k-\xi_{k+1}).
\end{equation*}
Direct use of Lemma \ref{bc} gives
\begin{equation*}
    \xi_k = O(\frac{1}{k}).
\end{equation*}
Using the fact $F(x^k)-\min F\leq \xi_k$, we then prove the result.

\section{Proof of Lemma \ref{lem-b}}
First, we prove that $x-\gamma_k\nabla f(x)$ is a contractive operator. For any $x,y$,
\begin{align*}
    &\|x-\gamma_k\nabla f(x)-y+\gamma_k\nabla f(y)\|^2\nonumber\\
    &\quad=\|x-y\|^2-2\gamma_k\langle\nabla f(x)-\nabla f(y),x-y\rangle\nonumber\\
    &\quad+\gamma_k^2\|\nabla f(x)-\nabla f(y)\|^2\nonumber\\
    &\quad\leq\|x-y\|^2-(\frac{2\gamma_k}{L}-\gamma_k^2)\|\nabla f(x)-\nabla f(y)\|^2\nonumber\\
    &\quad\leq\|x-y\|^2,
\end{align*}
where the first inequality depends on the fact $\langle\nabla f(x)-\nabla f(y),x-y\rangle\geq \frac{1}{L}\|\nabla f(x)-\nabla f(y)\|^2$, and the second one is due to $0<\gamma_k\leq\frac{2}{L}$.
Let $x^*$ be a minimizer of $F$. Obviously, it holds
\begin{equation*}
    x^*=\textbf{prox}_{\gamma_k g}[x^*-\gamma_k \nabla f(x^*)].
\end{equation*}
Noting $\textbf{prox}_{\gamma_k g}(\cdot)$ is  contractive,
\begin{align*}
&\|x^{k+1}-x^*\|=\|\textbf{prox}_{\gamma_k g}[x^k-\gamma_k \nabla f(x^k)\nonumber\\
&\quad+\beta_k(x^k-x^{k-1})]-\textbf{prox}_{\gamma_k g}[x^*-\gamma_k \nabla f(x^*)]\|\nonumber\\
&\quad\leq \|[x^k-\gamma_k \nabla f(x^k)+\beta_k(x^k-x^{k-1})]-[x^*-\gamma_k \nabla f(x^*)]\|\nonumber\\
&\quad\leq\|[x^k-\gamma_k \nabla f(x^k)]-[x^*-\gamma_k \nabla f(x^*)]\|\nonumber\\
&\quad+\|\beta_k(x^k-x^*+x^*-x^{k-1})\|\nonumber\\
&\quad\leq \|x^k-x^*\|+\beta_k\|x^k-x^*\|+\beta_k\|x^{k-1}-x^*\|.
\end{align*}
With Lemma \ref{app2}, we then prove the result.

\section{Proof of Theorem \ref{th-b}}
With Lemma \ref{lem-b}, the sequence is bounded. And it is easy to verify the boundedness of $\varepsilon_k$.  Thus, $\left\{\varepsilon_k(2\|x^{k+1}-\overline{x^{k+1}}\|^2+\|x^{k-1}-x^{k}\|^2)\right\}_{k\geq 1}$ is bounded. With almost the same proofs of Theorem \ref{sub-n}, we then prove the result.

\section{Proof of Theorem \ref{thm-linear}}
With (\ref{condtion}), we have
\begin{equation*}
    2\|x^{k+1}-\overline{x^{k+1}}\|^2\leq \frac{2}{\nu}( F(x^{k+1})-\min F)\leq \frac{2}{\nu}\xi_{k+1}\leq \frac{2}{\nu}\xi_k.
\end{equation*}
On the other hand, from the definition of (\ref{Lyapunov}),
\begin{equation*}
    \|x^k-x^{k-1}\|^2\leq  \frac{1}{\delta_k}  \xi_{k}.
\end{equation*}
With Lemma \ref{lem-sub}, we then derive
\begin{equation*}
     \xi_{k+1}^2\leq \varepsilon_k(\frac{1}{\delta_k}+\frac{2}{\nu})(\xi_k-\xi_{k+1})\cdot\xi_k.
\end{equation*}
With the assumption, $\sup_{k}\left\{\{\varepsilon_k(\frac{1}{\delta_k}+\frac{2}{\nu})\}_{k\geq 0}\right\}<+\infty$, and the bound is  assumed as $\ell>0$. And then, we have the following result,
\begin{equation*}
    \xi_{k+1}^2\leq \ell(\xi_k-\xi_{k+1})\cdot\xi_k.
\end{equation*}
If $ \xi_{k}=0$, we have $0= \xi_{k+1}=\xi_{k+2}=\ldots$. The result certainly holds. If $ \xi_{k}\neq0$,
\begin{equation*}
    (\frac{\xi_{k+1}}{\xi_{k}})^2+\ell(\frac{\xi_{k+1}}{\xi_{k}})-\ell\leq 0.
\end{equation*}
With basic algebraic computation,
\begin{equation*}
    \frac{ \xi_{k+1}}{ \xi_{k}}\leq \frac{2\ell}{\sqrt{\ell^2+4\ell}+\ell}.
\end{equation*}
By defining $\omega=\frac{2\ell}{\sqrt{\ell^2+4\ell}+\ell}$, we then prove the result.

\section{Proof of Lemma \ref{md-le1}}
For any $i\in[1,2,\ldots,m]$,
\begin{equation}\label{md-le1-t1}
    \frac{x^k_i-x^{k+1}_i}{\gamma_{k,i}}-\nabla_i^k f+\frac{\beta_{k,i}}{\gamma_{k,i}}(x^k_i-x^{k-1}_i)\in \partial g_i(x^{k+1}_i).
\end{equation}
With the convexity of $g_i$, we have
\begin{align}\label{md-le1-t2}
    &g_i(x^{k+1}_i)-g_i(x^k_i)\nonumber\\
    &\leq \langle \frac{x^{k+1}_i-x^k_i}{\gamma_{k,i}}+\nabla_i^k f+\frac{\beta_{k,i}}{\gamma_{k,i}}(x^{k-1}_i-x^k_i), x^k_i-x^{k+1}_i\rangle.
\end{align}
With (\ref{lpm}), we can have
\begin{align}\label{md-le1-t3}
    &f(x_1^{k+1},\ldots,x_{i-1}^{k+1},x_i^{k+1},x_{i+1}^{k}\ldots,x_m^k)\nonumber\\
    &\quad\quad-f(x_1^{k+1},\ldots,x_{i-1}^{k+1},x_i^{k},x_{i+1}^{k}\ldots,x_m^k)\nonumber\\
    &\leq \langle -\nabla_i^k f,x^k_i-x^{k+1}_i\rangle+\frac{L}{2}\|x^{k+1}-x^k\|^2.
\end{align}
Combining (\ref{md-le1-t2}) and (\ref{md-le1-t3}),
\begin{align}\label{md-le1-t4}
    &\left[f(x_1^{k+1},\ldots,x_{i-1}^{k+1},x_i^{k+1},x_{i+1}^{k}\ldots,x_m^k)+g_i(x^{k+1}_i)\right]\nonumber\\
    &\quad\quad-\left[f(x_1^{k+1},\ldots,x_{i-1}^{k+1},x_i^{k},x_{i+1}^{k}\ldots,x_m^k)+g_i(x^{k}_i)\right]\nonumber\\
    &\overset{(\ref{md-le1-t2})+(\ref{md-le1-t3})}{\leq} \frac{\beta_{k,i}}{\gamma_{k,i}}\langle x^{k}_i-x^{k-1}_i, x^{k+1}_i-x^k_i\rangle\nonumber\\
    &\quad\quad+(\frac{L}{2}-\frac{1}{\gamma_{k,i}})\|x^{k+1}_i-x^k_i\|^2\nonumber\\
    &\overset{a)}{\leq}\frac{\beta_{k,i}}{2\gamma_{k,i}}\|x^{k}_i-x^{k-1}_i\|^2+(\frac{L_i}{2}-\frac{1}{\gamma_{k,i}}+\frac{\beta_{k,i}}{2\gamma_{k,i}})\|x^{k+1}_i-x^k_i\|^2.
 \end{align}
where $a)$ uses the Schwarz inequality $\langle x^{k}_i-x^{k-1}_i, x^{k+1}_i-x^k_i\rangle\leq \frac{1}{2}\|x^{k}_i-x^{k-1}_i\|^2+\frac{1}{2}\|x^{k+1}_i-x^k_i\|^2$.
Summing (\ref{md-le1-t4}) from $i=1$ to $m$ yields
\begin{align}\label{md-le1-t5}
&F(x^{k+1})-F(x^k)\leq\sum_{i=1}^m\frac{\beta_{k,i}}{2\gamma_{k,i}}\|x^{k}_i-x^{k-1}_i\|^2\nonumber\\
&\quad\quad+\sum_{i=1}^m(\frac{L_i}{2}-\frac{1}{\gamma_{k,i}}+\frac{\beta_{k,i}}{2\gamma_{k,i}})\|x^{k+1}_i-x^k_i\|^2.
 \end{align}
With direct calculations and the non-increasity of $(\beta_{k,i})_{k\geq 0}$, we then obtain (\ref{nresult}).

\section{Proof of Lemma \ref{lem-sub-m}}
With  Lemma \ref{md-le1}, direct computing yields
\begin{align}\label{sketch-1-m}
   & \hat{\xi}_k- \hat{\xi}_{k+1}\geq\sum_{i=1}^m \frac{1}{2}(\frac{1-\beta_{k,i}}{\gamma_{k,i}}-\frac{L_i}{2})\nonumber\\
    &\quad\times(\|x^{k+1}_i-x^k_i\|^2+\|x^k_i-x^{k-1}_i\|^2)\nonumber\\
    &\quad=\frac{\underline{L}}{4}(\frac{1}{c}-1)\times(\|x^{k+1}-x^k\|^2+\|x^k-x^{k-1}\|^2).
\end{align}
For any $i\in[1,2,\ldots,m]$, the convexity of $g_i$ gives
\begin{align}\label{lem-sub-t3-m}
    g_i(x^{k+1}_i)- g_i([\overline{x^{k+1}}]_i)\leq\langle\widetilde{\nabla} g_i(x^{k+1}_i),x^{k+1}_i-[\overline{x^{k+1}}]_i\rangle,
\end{align}
where $\widetilde{\nabla} g_i(x^{k+1}_i)\in \partial g_i(x^{k+1}_i)$ and $[\overline{x^{k}}]_i$ denotes the $i$th coordinate of $\overline{x^{k}}$. With (\ref{md-le1-t1}), we then have
\begin{align}\label{lem-sub-t4-m}
    &g_i(x^{k+1}_i)- g_i([\overline{x^{k+1}}]_i)\nonumber\\
    &\leq\langle     \frac{x_i^k-x_i^{k+1}}{\gamma_{k,i}}-\nabla_i^k f+\frac{\beta_{k,i}}{\gamma_{k,i}}(x^k_i-x^{k-1}_i), x^{k+1}_i-[\overline{x^{k+1}}]_i\rangle.
\end{align}
Summing (\ref{lem-sub-t4-m}) with respect to $i$ from $1$ to $m$, and
\begin{equation}\label{lem-sub-t5-m}
    f(x^{k+1})-f(\overline{x^{k+1}})\leq \langle\nabla f(x^{k+1}),x^{k+1}-\overline{x^{k+1}}\rangle,
\end{equation}
we then have
\begin{align}\label{lem-sub-t6-m}
    &F(x^{k+1})-F(\overline{x^{k+1}})\nonumber\\
    &\quad\leq \sum_{i=1}^m\frac{\beta_{k,i}}{\gamma_{k,i}}\langle x^{k}_i-x^{k-1}_i,x^{k+1}_i-[\overline{x^{k+1}}]_i\rangle\nonumber\\
    &+\sum_{i=1}^m\langle  \frac{x^k_i-x^{k+1}_i}{\gamma_{k,i}},x^{k+1}_i-[\overline{x^{k+1}}]_i\rangle\nonumber\\
    &\quad+\sum_{i=1}^m\langle \nabla_i f(x^{k+1})-\nabla_i^k f,x^{k+1}_i-[\overline{x^{k+1}}]_i\rangle\nonumber\\
    &\overset{a)}{\leq}\sum_{i=1}^m\frac{\beta_{k,i}}{\gamma_{k,i}}\|x^{k}_i-x^{k-1}_i\|\cdot\|x^{k+1}_i-[\overline{x^{k+1}}]_i\|\nonumber\\
    &\quad+\sum_{i=1}^m \frac{1}{\gamma_{k,i}}\|  x^k_i-x^{k+1}_i\|\cdot\|x^{k+1}_i-[\overline{x^{k+1}}]_i\|\nonumber\\
    &+\sum_{i=1}^m L_i\|x^{k+1}-x^k\|\cdot\|x^{k+1}_i-[\overline{x^{k+1}}]_i\|\nonumber\\
    &\quad\overset{b)}{\leq}\sum_{i=1}^m\left(\frac{\|x^k_i-x^{k+1}_i\|}{\gamma_{k,i}}+\frac{\|x^{k}_i-x^{k-1}_i\|}{\gamma_{k,i}}+L_i\|x^{k+1}-x^k\|\right)\nonumber\\
    &\times\|x^{k+1}_i-[\overline{x^{k+1}}]_i\|,
\end{align}
where $a)$ is due to the Schwarz inequalities and the smooth assumption \textbf{A1}, $b)$ depends on the fact $0\leq\beta_{k,i}<1$. With (\ref{Lyapunov}) and (\ref{lem-sub-t6}), we have
\begin{align}\label{lem-sub-t7-m}
    &\hat{\xi}_{k+1}\leq\sum_{i=1}^m\left(\frac{\|x^{k+1}_i-x^{k}_i\|}{\gamma_{k,i}}+\frac{\|x^{k}_i-x^{k-1}_i\|}{\gamma_{k,i}}+L_i\|x^{k+1}-x^k\|\right)\nonumber\\
    &\quad\times\|x^{k+1}_i-[\overline{x^{k+1}}]_i\|+\sum_{i=1}^m \delta_{k+1,i}\|x^{k+1}_i-x^{k}_i\|^2.
\end{align}
Let
 \begin{align*}
    \hat{a}^k:=\begin{bmatrix}
                              \frac{1}{\gamma_{k,1}}\|x^{k+1}_1-x^k_1\| \\
                              \vdots\\
                              \frac{1}{\gamma_{k,m}}\|x^{k+1}_m-x^k_m\| \\
                              \frac{1}{\gamma_{k,1}}\|x^{k}_1-x^{k-1}_1\| \\
                              \vdots \\
                              \frac{1}{\gamma_{k,m}}\|x^{k}_m-x^{k-1}_m\| \\
                              L_1\|x^{k+1}-x^k\|\\
                              \vdots\\
                              L_m\|x^{k+1}-x^k\|\\
                               \delta_{k+1,1}\|x^{k+1}-x^{k}\| \\
                               \vdots \\
                                \delta_{k+1,m}\|x^{k+1}-x^{k}\| \\
\end{bmatrix},  \hat{b}^k:=\begin{bmatrix}
                              \|x^{k+1}_1-[\overline{x^{k+1}}]_1\| \\
                              \vdots\\
                               \|x^{k+1}_m-[\overline{x^{k+1}}]_m\| \\
                                 \|x^{k+1}_1-[\overline{x^{k+1}}]_1\| \\
                              \vdots\\
                               \|x^{k+1}_m-[\overline{x^{k+1}}]_m\| \\
                                 \|x^{k+1}_1-[\overline{x^{k+1}}]_1\| \\
                              \vdots\\
                               \|x^{k+1}_m-[\overline{x^{k+1}}]_m\| \\
                               \|x^{k+1}_1-x^{k}_1\| \\
                               \vdots\\
                               \|x^{k+1}_m-x^{k}_m\| \\
\end{bmatrix}.
\end{align*}
Using this and the definition of $\xi_{k+1}$ (\ref{Lyapunov}), we have:
\begin{align*}
    (\xi_{k+1})^2=\left|\langle \hat{a}^k,\hat{b}^k\rangle\right|^2\leq\|a^k\|^2\cdot\|b^k\|^2.
\end{align*}
Direct calculation yields
\begin{align*}
    &\|\hat{a}^k\|^2\leq\max\{\sum_{i=1}^m\left(\delta_{k+1,i}^2+L_i^2\right),\sum_{i=1}^m\frac{1}{\gamma^2_{k,i}}\}\nonumber\\
    &\quad\times(\|x^k-x^{k+1}\|^2+\|x^{k-1}-x^{k}\|^2)
\end{align*}
and
\begin{equation*}
    \|\hat{b}^k\|^2\leq 3\|x^{k+1}-\overline{x^{k+1}}\|^2+\|x^{k+1}-x^{k}\|^2.
\end{equation*}
Thus, we derive
\begin{align}\label{sketch-2-m}
    (\hat{\xi}_{k+1})^2&\leq\max\{\sum_{i=1}^m\left(\delta_{k+1,i}^2+L_i^2\right),\sum_{i=1}^m\frac{1}{\gamma^2_{k,i}}\}\nonumber\\
    &\times(\|x^k-x^{k+1}\|^2+\|x^{k-1}-x^{k}\|^2)\nonumber\\
    &\times(3\|x^{k+1}-\overline{x^{k+1}}\|^2+\|x^{k+1}-x^{k}\|^2).
\end{align}
Combining (\ref{sketch-1-m}) and (\ref{sketch-2-m}), we then prove the result.

\section{Proof of Theorem \ref{them-sub-m}}
With Lemma \ref{lem-sub-m},  $\sup_{k}\{\hat{\xi}_{k}\}<+\infty$, thus, $\sup_{k}\{F(x^k)\}<+\infty$ and $\sup_{k}\{\|x^{k}-x^{k-1}\|^2\}<+\infty$. Noting the coercivity of $D$, sequences $\{x^k\}_{k\geq 0}$ and $\{\overline{x^k}\}_{k\geq 0}$ are  bounded. With the assumptions on $\gamma_{k,i}$ and $\beta_{k,i}$, $\sup_k\{\hat{\varepsilon}_k\}<+\infty$. Thus, $\left\{\hat{\varepsilon}_k(3\|x^{k+1}-\overline{x^{k+1}}\|^2+\|x^{k-1}-x^{k}\|^2)\right\}_{k\geq 1}$ is bounded; and we assume the bound is $R$, i.e.,
\begin{equation*}
    \sup_{k}\left\{\hat{\varepsilon}_k(3\|x^{k+1}-\overline{x^{k+1}}\|^2+\|x^{k-1}-x^{k}\|^2)\right\}\leq R.
\end{equation*}
With Lemma \ref{lem-sub-m},  we then have
\begin{equation*}
    \hat{\xi}_{k+1}^2\leq R(\hat{\xi}_k-\hat{\xi}_{k+1}).
\end{equation*}
From Lemma \ref{bc},
\begin{equation*}
    \hat{\xi}_k= O(\frac{1}{k}).
\end{equation*}
Using the fact $F(x^k)-\min F\leq \hat{\xi}_k$, we then obtain the result.

\section{Proof of Theorem \ref{them-linear-m}}
With the optimal strong
convexity condition, we have
\begin{equation*}
    3\|x^{k+1}-\overline{x^{k+1}}\|^2\leq \frac{3}{\nu}\hat{\xi}_{k+1}\leq\frac{3}{\nu}\hat{\xi}_{k}.
\end{equation*}
The direct computing yields
\begin{equation*}
    \|x^{k-1}-x^{k}\|^2\leq\frac{\hat{\xi}_{k}}{\min_i\{\delta_{k,i}\}}
\end{equation*}
With Lemma \ref{lem-sub-m},
\begin{equation*}
    (\hat{\xi}_{k+1})^2\leq\left(\hat{\varepsilon}_k+\frac{3}{\nu}+\frac{1}{\min_i\{\delta_{k,i}\}}\right)(\hat{\xi}_k-\hat{\xi}_{k+1})\hat{\xi}_{k}.
\end{equation*}
It is easy to see that $\hat{\varepsilon}_k+\frac{3}{\nu}+\frac{1}{\min_i\{\delta_{k,i}\}}$ is bounded by some positive constant $\ell>0$. Thus we have
\begin{equation*}
    F(x^k)-\min F\leq \hat{\xi}_k= O\left((\frac{2\ell}{\sqrt{\ell^2+4\ell}+\ell})^k\right).
\end{equation*}
Letting $\omega=\frac{2\ell}{\sqrt{\ell^2+4\ell}+\ell}$, we then prove the result.

\section{Proof of Lemma \ref{ms-le1}}
In the $k$-th iteration, it holds
\begin{equation*}\label{ms-le1-t1}
    \frac{x^k_{i_k}-x^{k+1}_{i_k}}{\gamma_{k}}-\nabla_{i_k} f(x^k)+\frac{\beta_{k}}{\gamma_{k}}(x^k_{i_k}-x^{k-1}_{i_k})\in \partial g_{i_k}(x^{k+1}_{i_k}).
\end{equation*}
With the convexity of $g_{i_k}$, we have
\begin{align}\label{ms-le1-t2}
    &g_{i_k}(x^{k+1}_{i_k})-g_{i_k}(x^k_{i_k})\nonumber\\
    &\leq \langle \frac{x^{k+1}_{i_k}-x^k_{i_k}}{\gamma_{k}}+\nabla_{i_k} f(x^k)+\frac{\beta_{k}}{\gamma_{k}}(x^{k-1}_{i_k}-x^k_{i_k}), x^k_{i_k}-x^{k+1}_{i_k}\rangle
\end{align}
With the Lipschitz of $\nabla f$, we then derive
\begin{align}\label{ms-le1-t3}
    &f(x^{k+1})-f(x^k)\nonumber\\
    &\quad\leq \langle -\nabla_{i_k} f(x^k), x^k_{i_k}-x^{k+1}_{i_k}\rangle+\frac{L}{2}\|x^{k+1}-x^k\|^2,
\end{align}
where we used the fact $\langle -\nabla_{i_k} f(x^k),x^k_{i_k}-x^{k+1}_{i_k}\rangle=\langle -\nabla f(x^k),x^k-x^{k+1}\rangle$.
Combining (\ref{ms-le1-t2}) and (\ref{ms-le1-t3}),
\begin{align}\label{ms-le1-t4}
    &F(x^{k+1})-F(x^k)\nonumber\\
    &\quad=\left[f(x^{k+1})+g_{i_{k}}(x^{k+1}_{i_k})\right]-\left[f(x^{k})+g_{i_{k}}(x^{k}_{i_k})\right]\nonumber\\
    &\overset{a)}{\leq} \frac{\beta_{k}}{\gamma_{k}}\langle x^{k}_{i_k}-x^{k-1}_{i_k}, x^{k+1}_{i_k}-x^k_{i_k}\rangle+(\frac{L}{2}-\frac{1}{\gamma_{k}})\|x^{k+1}-x^k\|^2\nonumber\\
    &\quad\overset{b)}{\leq}(\frac{L}{2}-\frac{1}{\gamma_{k}}+\frac{\beta_{k}}{2\sqrt{m}\gamma_{k}})\|x^{k+1}-x^k\|^2\nonumber\\
    &+\frac{\sqrt{m}\beta_{k}}{2\gamma_{k}}\|x^{k}_{i_k}-x^{k-1}_{i_k}\|^2.
 \end{align}
where $a)$ depends on (\ref{ms-le1-t2})+(\ref{ms-le1-t3}), and $b)$ uses the Schwarz inequality $\langle x^{k}_i-x^{k-1}_i, x^{k+1}_i-x^k_i\rangle\leq \frac{\sqrt{m}}{2}\|x^{k}_i-x^{k-1}_i\|^2+\frac{1}{2\sqrt{m}}\|x^{k+1}_i-x^k_i\|^2$ and the fact $\|x^{k+1}-x^k\|^2=\|x^{k+1}_{i_k}-x^k_{i_k}\|^2$.
Taking conditional conditional expectations of (\ref{ms-le1-t4}) on $\chi^k$,
\begin{align}\label{ms-le1-t5}
&\EE[F(x^{k+1})\mid\chi^k]-F(x^k)\leq\frac{\beta_{k}}{2\sqrt{m}\gamma_{k}}\|x^{k}-x^{k-1}\|^2\nonumber\\
&\quad\quad+(\frac{L}{2}-\frac{1}{\gamma_{k}}+\frac{\beta_{k}}{2\sqrt{m}\gamma_{k}})\EE(\|x^{k+1}
-x^k\|^2\mid\chi^k).
 \end{align}
Taking total expectations on (\ref{ms-le1-t5}), and using $\EE(\EE(\cdot\mid\chi^k))=\EE(\cdot)$,
\begin{align*}
&\EE F(x^{k+1}) -\EE F(x^k)\leq\frac{\beta_{k}}{2\sqrt{m}\gamma_{k}}\EE\|x^{k}-x^{k-1}\|^2\\
&\quad+(\frac{L}{2}-\frac{1}{\gamma_{k}}+\frac{\beta_{k}}{2\sqrt{m}\gamma_{k}})\EE\|x^{k+1}
-x^k\|^2.
 \end{align*}
 Thus  we have
 \begin{align*}
    & \left[\EE F(x^{k})+\frac{\beta_k}{2\sqrt{m}\gamma_k}\EE\|x^{k}-x^{k-1}\|^2\right]\\
    &\quad\quad-\left[\EE F(x^{k+1})+\frac{\beta_{k}}{2\sqrt{m}\gamma_{k}}\EE\|x^{k+1}-x^{k}\|^2\right]\\
     &\geq(\frac{1-\beta_k/\sqrt{m}}{\gamma_k}-\frac{L}{2})\EE\|x^{k+1}-x^k\|^2.
\end{align*}
With the non-increasity of $\{\beta_k\}_{k\geq 0}$, $\left\{\frac{\beta_k}{2\sqrt{m}\gamma_k}=\frac{\beta_k L}{4(1-\beta_k/\sqrt{m})c}\right\}_{k\geq 0}$ is also non-increasing; and then, we get
\begin{align}\label{litpre}
\sum_{k}\EE\|x^{k+1}-x^k\|^2<+\infty.
\end{align}
 Using \eqref{litpre} and [Lemma 3; \cite{davis2016convergence}], we can obtain
\begin{align}\label{lit-yin}
 \min_{0\leq i\leq k}\EE(\|x^{i+1}-x^i\|^2)=o(\frac{1}{k}).
\end{align}
Noticing that
$
\min_{0\leq i\leq k}\|x^{i+1}-x^i\|^2\leq \|x^{i+1}-x^i\|^2,~~i\in\{0,1,2\ldots,k\}
$, we are then led to
\begin{align*}
\EE(\min_{0\leq i\leq k}\|x^{i+1}-x^i\|^2)\leq \EE\|x^{i+1}-x^i\|^2,~~i\in\{0,1,2\ldots,k\},
\end{align*}
which together with \eqref{lit-yin} tells us
\begin{align}
\EE(\min_{0\leq i\leq k}\|x^{i+1}-x^i\|^2)\leq \min_{0\leq i\leq k}\EE\|x^{i+1}-x^i\|^2=o(\frac{1}{k}).
\end{align}

\section{Proof of Theorem \ref{ms-th1}}
Direct calculations yield
\begin{align*}
&\|S_{\gamma}(x^k)\|^2\nonumber\\
&=m\cdot\EE(\|x^k_{i_k}-\textbf{prox}_{\gamma g_{i_k}}[x^k_{i_k}-\gamma\nabla_{i_k} f(x^k)]\|^2\mid\chi^k)\nonumber\\
&=m\cdot \EE(\|x^k_{i_k}-\textbf{prox}_{\gamma g_{i_k}}[x^k_{i_k}-\gamma\nabla_{i_k} f(x^k)+\beta(x^k_{i_k}-x^{k-1}_{i_k})]\nonumber\\
&+\textbf{prox}_{\gamma g_{i_k}}[x^k_{i_k}-\gamma\nabla_{i_k} f(x^k)+\beta(x^k_{i_k}-x^{k-1}_{i_k})]\nonumber\\
&-\textbf{prox}_{\gamma g_{i_k}}[x^k_{i_k}-\gamma\nabla_{i_k} f(x^k)]\|^2\mid\chi^k)\nonumber\\
&\leq 2m\cdot\EE(\|x^{k+1}-x^k\|^2\mid\chi^k)+2m\cdot\EE(\|x^{k}_{i_k}-x^{k-1}_{i_k}\|^2\mid\chi^k)\nonumber\\
&=2m\cdot\EE(\|x^{k+1}-x^k\|^2\mid\chi^k)+2\|x^{k}-x^{k-1}\|^2.\nonumber
\end{align*}
Taking expectations of both sides,
\begin{align*}
\EE\|S_{\gamma}(x^k)\|^2\leq 2m\EE\|x^{k+1}-x^k\|^2+2\EE\|x^{k}-x^{k-1}\|^2.
\end{align*}
From \eqref{ms-result-2}, the result is then  proved.

\section{Proof of Lemma \ref{lem:linear-convergence}}

The optimization condition of  iteration \eqref{schemes} yields
\begin{align}\label{linear-convergence-t1}
   &g_{i_k}([\overline{x^k}]_{i_k})- g_{i_k}(x^{k+1}_{i_k})\nonumber\\
   &\quad+\langle (\overline{x^k}-x^{k+1})_{i_k},\nabla_{i_k} f(x^k)+\frac{\beta}{\gamma}(x^k_{i_k}-x^{k-1}_{i_k})\rangle\nonumber\\
   &\geq\frac{1}{\gamma}\langle (\overline{x^k}-x^{k+1})_{i_k},x^k_{i_k}-x^{k+1}_{i_k}\rangle.
\end{align}
With the scheme of the algorithm, we then get
\begin{align*}
    &\|x^{k+1}-\overline{x^{k+1}}\|^2 \leq \|x^{k+1}-\overline{x^{k}}\|^2\nonumber\\
    &\qquad=\|x^k-\overline{x^{k}}\|^2+\|x^{k+1}-x^k\|^2\\
    &-2\langle (x^k-\overline{x^{k}})_{i_k}, (x^{k}-x^{k+1})_{i_k}\rangle\nonumber\\
    &\qquad=\|x^k-\overline{x^{k}}\|^2-\|x^{k+1}-x^k\|^2\\
    &+2\langle (\overline{x^{k}}-x^{k+1})_{i_k}, (x^{k}-x^{k+1})_{i_k}\rangle\nonumber\\
   (\ref{linear-convergence-t1}) &\leq\|x^k-\overline{x^{k}}\|^2-\|x^{k+1}-x^k\|^2
   &\qquad+2\beta\langle (\overline{x^{k}}-x^{k+1})_{i_k}, (x^{k}-x^{k-1})_{i_k}\rangle\nonumber\\
    &+2\gamma[\langle(\overline{x^{k}}-x^{k+1})_{i_k},\nabla_{i_k} f(x^k)\rangle\\
    &\qquad+g_{i_k}([\overline{x^{k}}]_{i_k})-g_{i_k}(x^{k+1}_{i_k})]\nonumber\\
    &=\|x^k-\overline{x^{k}}\|^2-\|x^{k+1}-x^k\|^2\nonumber\\
    &\qquad+\underbrace{2\gamma\langle(\overline{x^{k}}-x^k)_{i_k},\nabla_{i_k}f(x^k)\rangle}_{\textrm{I}}\nonumber\\
    &+\underbrace{2\gamma[\langle(x^{k}-x^{k+1})_{i_k},\nabla_{i_k} f(x^k)\rangle+g_{i_k}([\overline{x^{k}}]_{i_k})-g_{i_k}(x^{k+1}_{i_k})]}_{\textrm{II}}\nonumber\\
    &+\underbrace{2\beta\langle (\overline{x^{k}}-x^{k+1})_{i_k}, (x^{k}-x^{k-1})_{i_k}\rangle}_{\textrm{III}}\nonumber
\end{align*}
In the following, we bound the expectations of I, II and III:
\begin{align*}
  &\textrm{I}=2\gamma\EE(\langle(\overline{x^{k}}-x^k)_{i_k},\nabla_{i_k}f(x^k)\rangle\mid\chi^k)\\
  &\qquad=\frac{2\gamma}{m}\langle\overline{x^{k}}-x^k,\nabla f(x^k)\rangle\leq\frac{2\gamma}{m}(f(\overline{x^{k}})-f(x^k)).
\end{align*}
Noticing that  $\langle(x^{k}-x^{k+1})_{i_k},\nabla_{i_k} f(x^k)\rangle=\langle x^{k}-x^{k+1},\nabla f(x^k)\rangle\leq f(x^k)-f(x^{k+1})+\frac{L}{2}\|x^{k+1}-x^k\|^2$, we can obtain
\begin{align*}
   \textrm{II}\leq&2\gamma\EE(f(x^k)-f(x^{k+1})+\frac{L}{2}\|x^{k+1}-x^k\|^2\nonumber\\
   &\qquad+g_{i_k}([\overline{x^{k}}]_{i_k})-g_{i_k}(x^{k+1}_{i_k})\mid\chi^k)\nonumber\\
    &=2\gamma f(x^k)-2\gamma\EE(f(x^{k+1})|\chi^k)+\frac{2\gamma}{m}g(\overline{x^{k}})\nonumber\\
    &\qquad-2\gamma\EE(g(x^{k+1})\mid\chi^k)\nonumber\\
    &+2\gamma\frac{m-1}{m}g(x^k)+\gamma L\EE\|x^{k+1}-x^k\|^2.\nonumber
\end{align*}
With Schwarz inequality,
\begin{align*}
  \textrm{III}=&2\beta \EE(\langle (\overline{x^{k}}-x^{k+1})_{i_k}, (x^{k}-x^{k-1})_{i_k}\rangle\mid\chi^k)\nonumber\\
  &\leq\beta\EE(\|\overline{x^{k}}-x^{k+1}\|^2\mid\chi^k)+\beta\EE(\|x^{k}-x^{k-1}\|^2\mid\chi^k)\nonumber\\
  &\leq2\beta\|\overline{x^{k}}-x^{k}\|^2+2\beta\EE(\|x^{k}-x^{k+1}\|^2\mid\chi^k)\nonumber\\
  &+\frac{\beta}{m}\|x^{k}-x^{k-1}\|^2.
\end{align*}
Combining the bounds of I, II and III,
\begin{align*}
&\EE(\|x^{k+1}-\overline{x^{k+1}}\|^2\mid\chi^k)\leq (1+2\beta)\|x^{k}-\overline{x^{k}}\|^2\\
&\quad\quad-(1-\gamma L-2\beta)\EE(\|x^{k+1}-x^k\|^2\mid\chi^k)\\
&+\frac{2\gamma}{m}(F^*-F(x^k))+2\gamma(F(x^k)-\EE(F(x^{k+1})\mid\chi^k))\\
&\quad\quad+\frac{\beta}{m}\|x^{k}-x^{k-1}\|^2.
\end{align*}
Taking expectations gives
\begin{align}\label{linear-convergence-t2}
&\frac{2\gamma}{m}\EE(F(x^k)-F^*)-2\beta\EE\|x^{k}-\overline{x^{k}}\|^2\nonumber\\
&\quad\leq \EE[\|x^{k}-\overline{x^{k}}\|_2^2+2\gamma F(x^{k})]\nonumber\\
&-\EE[\|x^{k+1}-\overline{x^{k+1}}\|_2^2+2\gamma F(x^{k+1})]\nonumber\\
&\quad-(1-\gamma L-2\beta)\EE\|x^{k+1}-x^k\|_2^2+\frac{\beta}{m}\EE\|x^{k}-x^{k-1}\|^2.
\end{align}
With condition \eqref{condtion}, we can further obtain
\begin{align}\label{linear-convergence-t3}
&\frac{2\gamma}{m}(F(x^k)-F^*)-2\beta\EE\|x^{k}-\overline{x^{k}}\|^2\nonumber\\
&\quad=\frac{\gamma}{m}(F(x^k)-F^*)+\frac{\gamma}{m}(F(x^k)-F^*)-2\beta\|x^{k}-\overline{x^{k}}\|^2\nonumber\\
&=\frac{\gamma}{m}(F(x^k)-F^*)+(\frac{\gamma\nu}{m}-2\beta)\|x^{k}-\overline{x^{k}}\|^2\nonumber\\
&\quad\geq\underline{\ell}\cdot\EE\left[\frac{\gamma}{m}(F(x^{k})-F^*)+\|x^{k}-\overline{x^{k}}\|^2\right].
\end{align}
With (\ref{linear-convergence-t2}) and (\ref{linear-convergence-t3}), we derive the result.

\section{Proof of Theorem \ref{th:linear-convergence}}
Direct computation gives
\begin{equation*}
    \gamma L+2\beta-1\leq(1-\underline{\ell})\frac{\beta}{m}~,\,~0<\underline{\ell}<1.
\end{equation*}
Thus, from Lemma \ref{lem:linear-convergence},
\begin{align*}
&\underline{\ell}\EE[\|x^{k}-\overline{x^{k}}\|^2+2\gamma(F(x^{k})-F^*)+\frac{\beta}{m}\EE\|x^{k}-x^{k-1}\|^2] \\
&\quad\leq\EE[\|x^{k}-\overline{x^{k}}\|^2+2\gamma(F(x^{k})-F^*)+\frac{\beta}{m}\|x^{k}-x^{k-1}\|^2]\nonumber\\
&-\EE[\|x^{k+1}-\overline{x^{k+1}}\|^2+2\gamma(F(x^{k+1})-F^*)\\
&\quad+\frac{\beta}{m}\|x^{k+1}-x^{k}\|^2].\nonumber
\end{align*}
That is also
\begin{align*}
&\EE\Big[\|x^{k+1}-\overline{x^{k+1}}\|^2+2\gamma(F(x^{k+1})-F^*)\\
&+\frac{\beta}{m}\|x^{k+1}-x^{k}\|^2\Big]\leq(1-\underline{\ell})\EE\Big[\|x^{k}-\overline{x^{k}}\|^2\\
&+2\gamma(F(x^{k})-F^*)+\frac{\beta}{m}\EE\|x^{k}-x^{k-1}\|^2\Big].
\end{align*}
We are then led to
\begin{align*}
&\EE\Big[\|x^{k+1}-\overline{x^{k+1}}\|^2+2\gamma(F(x^{k+1})-F^*)\\
&\qquad+\frac{\beta}{m}\|x^{k+1}-x^{k}\|^2\Big]\leq(1-\underline{\ell})^k\times \EE\Big[\|x^{1}-\overline{x^{1}}\|^2\\
&+2\gamma(F(x^{1})-F^*)+\frac{\beta}{m}\EE\|x^{1}-x^{0}\|^2\Big],
\end{align*}
which yields
\begin{align*}
&\EE(2\gamma(F(x^{k+1})-F^*))\leq\EE\Big[\|x^{k+1}-\overline{x^{k+1}}\|^2\\
&+2\gamma(F(x^{k+1})-F^*)+\frac{\beta}{m}\|x^{k+1}-x^{k}\|^2\Big]=O\left((1-\frac{\gamma\nu}{2m})^k\right).
\end{align*}

\end{document}